\documentclass[11pt,leqno]{article}

\usepackage{amsmath,amsfonts,amscd,amssymb,theorem}

\long\def\comment#1\endcomment{}

\makeatletter
\begingroup
\gdef\th@dotted{\normalfont\itshape
  \def\@begintheorem##1##2{%
        \item[\hskip\labelsep \theorem@headerfont ##1\ ##2.]}%
\def\@opargbegintheorem##1##2##3{%
   \item[\hskip\labelsep \theorem@headerfont ##1\ ##2\ (##3).]}}
\endgroup
\makeatother

\theoremstyle{dotted}

\newtheorem{theorem}{Theorem}[section]
\newtheorem{lemma}[theorem]{Lemma}
\newtheorem{conj}[theorem]{Conjecture}

\newtheorem{corr}[theorem]{Corollary}

\makeatletter
\begingroup
\gdef\th@upshape{\normalfont
  \def\@begintheorem##1##2{%
        \item[\hskip\labelsep \theorem@headerfont ##1\ ##2.]}%
\def\@opargbegintheorem##1##2##3{%
   \item[\hskip\labelsep \theorem@headerfont ##1\ ##2\ (##3).]}}
\endgroup
\makeatother

\theoremstyle{upshape}

\newtheorem{remark}[theorem]{Remark}

\makeatletter
\renewcommand{\subsection}{\@startsection{subsection}{2}{0pt}{-3ex
plus -1ex minus -0.2ex}{-2mm plus -0pt minus
-2pt}{\normalfont\bfseries}} \makeatother

\makeatletter
\@addtoreset{equation}{section}
\makeatother

\newcommand{\cntrct}                
{\hspace{2pt}\raisebox{1pt}{\text{$\lrcorner$}}\hspace{2pt}}

\newcommand{\proof}[1][Proof.]{\smallskip\noindent{\em #1}}
\def\endproof{\hfill\ensuremath{\square}\par\medskip}

\def\eqref#1{\thetag{\ref{#1}}}

\let\latexref=\ref
\def\ref#1{{\normalfont{\latexref{#1}}}}

\newcommand{\wt}{\widetilde}

\newcommand{\idot}{{\:\raisebox{1pt}{\text{\circle*{1.5}}}}}
\newcommand{\hdot}{{\:\raisebox{3pt}{\text{\circle*{1.5}}}}}
\newcommand{\calo}{{\cal O}}

\newcommand{\gr}{\operatorname{{\sf gr}}}
\newcommand{\id}{{\operatorname{{\sf id}}}}
\newcommand{\age}{\operatorname{{\sf age}}}
\renewcommand{\dim}{\operatorname{{\sf dim}}}
\newcommand{\codim}{\operatorname{{\sf codim}}}
\newcommand{\graph}{\operatorname{{\sf graph}}}
\newcommand{\bimod}{{\operatorname{\sf bimod}}}

\newcommand{\C}{{\mathbb C}}
\newcommand{\Q}{{\mathbb Q}}
\newcommand{\Z}{{\mathbb Z}}

\newcommand{\A}{{\mathbb A}}

\newcommand{\g}{{\mathfrak g}}
\newcommand{\h}{{\mathfrak h}}

\newcommand{\ZZ}{{\bf\sf Z}}

\newcommand{\D}{{\mathcal D}}
\newcommand{\T}{{\mathcal T}}

\newcommand{\Ext}{\operatorname{Ext}}

\renewcommand{\Im}{\operatorname{{\sf Im}}}
\newcommand{\Hilb}{\operatorname{{\sf Hilb}}}

\newcommand{\Spec}{\operatorname{Spec}}

\title{Multiplicative McKay correspondence in the symplectic case}

\author{D. Kaledin\thanks{Partially supported by CRDF grant
RM1-2354-MO02.}}

\begin{document}

\maketitle

\tableofcontents

\section*{Introduction}

McKay correspondence as a subject appeared approximately twenty-five
years ago; for an excellent overview we refer the reader to
M. Reid's Bourbaki talk \cite{R2}. The basic setup is the
following. Let $V$ be a finite-dimensional complex vector space, and
let $G \subset SL(V)$ be a finite subgroup which acts on $V$
preserving the volume. The quotient $Y = V/G$ is a singular
algebraic variety with trivial canonical bundle. Assume given a
smooth projective resolution $X \to Y$ which is {\em crepant} --- in
other words, such that the canonical bundle $K_X$ is trivial. What
can one say about the topological invariants of $X$, starting from
the combinatorics of the $G$-action on $V$? Somewhat surprisingly,
it turns out that one can actually say quite a lot. This is, roughly
speaking, what is called the McKay correspondence.

For some time, all the research concentrated on the case of $\dim V
= 2$ (in particular, this was the case with the original paper
\cite{McK} of J. McKay which gave name to the whole subject).  An
important breakthrough was achieved in \cite{R1}. Drawing on the
case of $\dim V = 2$, well-understood by that time, and also on some
physics papers and on previous work in $\dim 3$, M. Reid has
proposed several conjectures for vector spaces of arbitrary
dimension. Essentially, he predicted what the homology
$H_\idot(X,\Q)$, the cohomology $H^\hdot(X,\Q)$ and the $K$-theory
$K_0(X)$ should be, and gave some pointers to possible proofs.

There has been considerable progress since then, and it mostly
concentrated on two areas. Firstly, most of the conjectures in
\cite{R1} have been proved (\cite{B}, \cite{DL}) by applying
M. Kontsevich's method of {\em motivic integration}. Secondly,
W. Chen and Y. Ruan \cite{CR} have introduced a new invariant for
the pair $\langle V, G \subset SL(V) \rangle$ called the {\em
orbifold cohomology algebra} $H_{orb}([V/G])$. As a vector space,
this ring coincides with what McKay correspondence predict must be
$H^\hdot(X)$.  However, and this is the new ingredient, it also
carries a structure of a commutative algebra.  A natural question,
then, would be to take a crepant resolution $X \to V/G$ and to
compare the multiplication $H^\hdot_{orb}([V/G])$ with the standard
multiplication in $H^\hdot(X)$.

Unfortunately, the only general technique in McKay correspondence,
--- that of motivic integration, --- does not seem to apply to this
problem. And on the other hand, in the general case it is not even
expected that the two multiplicative structures in
$H^\hdot_{orb}([V/G]) \cong H^\hdot(X)$ are the same (see \cite{Ru2}
for precise conjectures). In order to get results, one has to
restrict generality in some way.

\bigskip

The particular case that we will consider in this paper is that of
{\em symplectic} vector spaces --- we will assume that the vector
space $V$ is equiped with a non-degenerate skew-symmetric form, and
that the finite group $G \subset Sp(V) \subset SL(V)$ preserves not
only the volume, but also the symplectic form.

\bigskip

In this case, there are several simplifications. Firstly, the
additive McKay correspondence can be established explicitly, without
recourse to the motivic integration techniques (\cite{K}). One also
gets a basis in cohomology given by classes of explicit algebraic
cycles. Secondly, the orbifold cohomology space and the orbifold
cohomology algebra are much easier to describe than in the general
case. However, the most important simplification is that the
multiplication in $H^\hdot(X,\Q)$ is known at least in one
particular example --- one takes $V = \C^{2n}$, one takes $G$ to be
the symmetric group on $n$ letters, and one takes $X \to V/G$ to be
the Hilbert scheme of subschemes in $\C^2$ of dimension $0$ and
length $n$. In fact, the multiplicative structure in $H^\hdot(X,\Q)$
has been discovered indepedently in \cite{LS} and \cite{V} (see also
\cite{wang}) without any reference to the general orbifold
cohomology construction of \cite{CR}. It is only later that it was
realized that the cohomology algebras are one and the same. Y. Ruan
\cite{Ru1} has conjectured that this should always be true in the
symplectic case, no matter what particular crepant resolution one
takes.

The goal of this paper is to present a result proved recently by the
author together with V. Ginzburg. The result, which is \cite[Theorem
1.2]{GK}, establishes the Ruan conjecture: we prove that in the
symplectic case, the cohomology $H^\hdot(X,\C)$ of {\em any} crepant
resolution $X \to V/G$ coincides with the Chen-Ruan orbifold
cohomology algebra $H^\hdot_{orb}([V/G],\C)$.

Unfortunately, our task is made quite difficult by the nature of the
proof. The proof is really quite roundabout, and it relies heavily
on earlier results obtained in other papers and/or by other
authors. The identification $H^\hdot_{orb}([V/G],\C) \cong
H^\hdot(X,\C)$ is split into a series of algebra isomorphisms most
of which have been known before. Only some isomorphisms in the
middle still need proving --- and the things isomorphic are really
quite far removed from either $H^\hdot_{orb}([V/G])$ or
$H^\hdot(X)$. The upshot is that there is very little in \cite{GK}
which has to do with actual multiplication in cohomology, be it
orbifold cohomology or the usual one. For that matter, more than
half of the paper is taken with hardcore deformation theory, and it
even has nothing to do with symplectic quotient singularities or
crepant resulutions. This makes \cite{GK} hard to read for people
not familiar with the field.

Faced with this difficulty, we have adopted a solution which is
somewhat radical -- we hope that it is justified, if not by
pedagogical reasons, then at least by its comparative novelty. There
is an alternative prove of \cite[Theorem 1.2]{GK} which is shorter,
more conceptual and more to the point. However, it needs some
intermediate results that are presently not known. Worse than that,
some of them might not even be true. Nevertheless, it is this fake
``proof'' that we will present in this paper. By this, we hope to
achieve two goals: firstly, we show the main ideas behind the real
proof without bothering with making them work, secondly, we motivate
the introduction of some techniques, such as Hochschild cohomology,
which {\em a priori} would seem quite foreign to the subject. Then
in the last section, we describe how the actual proof in \cite{GK}
works --- what additional ingredients one has to add to the brew so
that all the ``fake'' statements of earlier sections become
theorems.

Therefore this paper should very much be thought of as an exposition
--- very little, if anything, will be actually proved. We hope that
at least, it will make a passable read.

\subsection*{Acknowledgements.} This is a write-up of my talk at
the Atelier sur les structures alg\'ebriques et espaces de modules
which took place at the University of Montreal in July 2003. I would
like to thank the organizers for inviting me to this great
conference and for giving me the opportunity to present our
results. The hospitality of the University of Montreal was also of
highest standard and is gratefully acknowledged.

\section{Orbifold cohomology algebra}

We start with combinatorics. Let $V$ be a finite dimensional complex
vector space. Let $g:V \to V$ be an automorphism of $V$ of some
finite order $r$. Then $V$ can be decomposed into $g$-eigenspaces,
$V = V_1 \oplus \dots \oplus V_k$, where $g$ acts on $V_i$ by
multiplication by $\exp(2\pi\sqrt{-1} \cdot a_i/r)$ for some integer
$a_i$, $0 \leq a_i < r$. Denote by $m_i = \dim V_i$ the multiplicity
of the eigenvalue $\exp(2\pi\sqrt{-1} \cdot a_i/r)$. By the {\em
age} of the automorphism $g$ one understands the sum
$$
\age g = \sum_{1 \leq i \leq k}m_i\frac{a_i}{r}.
$$
The age is {\em a priori} a rational number. However, the
automorphism $g$ acts on $\det V$ by $\exp(2\pi\sqrt{-1}\age
g)$. Therefore if $g$ preserves the volume on $V$, its age is
actually an integer. It is obviously invariant under conjugation.

Assume given a finite subgroup $G \subset SL(V)$, and let
$\overline{G}$ denote the set of conjugacy classes of elements in
$G$. By the {\em orbifold cohomology space}
$H^\hdot_{orb}([V/G],\C)$ we will understand the graded vector space
defined by
$$
H^k_{orb}([V/G],\C) = \bigoplus_{g \in \overline{G},k = 2\age g}\C
\cdot g.
$$
In other words, $H^\hdot_{orb}([V/G],\C)$ is generated by conjugacy
classes $g \in \overline{G}$, each placed in degree $2\age g$, and
the odd part of $H^\hdot_{orb}([V/G],\C)$ is trivial. One can also
define a bigger space $H^\hdot_G(V,\C)$ by
$$
H^k_G(V,\C) = \bigoplus_{g \in G,k = 2\age g}\C \cdot g.
$$
In other words, $H^\hdot_{G}(V,\C)$ is just the vector space $\C[G]$
spanned by the group $G$ and graded by twice the age. The group $G$
acts on $H^\hdot_G(V,\C)$ by conjugation, and
$H^\hdot_{orb}([V/G],\C)$ is naturally identified with the subspace
of $G$-invariant vectors in $H^\hdot_G(V,\C)$.

To define a multiplication on $H^\hdot_{orb}([V/G],\C)$ -- we follow
the exposition of \cite{CR} given in \cite{FG} -- one defines a
$G$-invariant multiplication on the bigger space
$H^\hdot_G(V,\C)$. This is given in the basis numbered by $g \in G$
by
$$
g_1 \cdot g_2 = c(g_1,g_2)g_1g_2,
$$
where $c(g_1,g_2)=1$ if $\age g_1 + \age g_2 = \age (g_1g_2)$ and
$0$ otherwise. It is extremely non-trivial to prove that this really
defines an associative multiplication on $H^\hdot_G(V,\C)$. Once one
has proved this, one descends to the subspace
$H^\hdot_{orb}([V/G],\C)$ of $G$-invariant vectors, and observes that
the multiplication becomes not only associative, but also
commutative.

Assume now that the vector space $V$ is equipped with a
non-degenerate skew-symmetric form $\omega \in \Lambda^2V^*$ and
that $G \subset Sp(V)$ preserves $\omega$. Then the orbifold
cohomology algebra is much easier to describe.

\begin{lemma}\label{filtr}
For every element $g \in G \subset Sp(V)$, we have
$$
\age g = \frac{1}{2}\codim V^g,
$$
where $V^g \subset V$ is the subspace of $g$-invariant
vectors. Moreover, define an increasing filtration $F_\idot\C[G]$ on
the group algebra $\C[G]$ by
$$
F_k\C[G] = \bigoplus_{g \in G,\codim V^g \leq k}\C \cdot g.
$$
Then this filtration is compatible with the group algebra
multiplication, and we have
$$
H^\hdot_G(V,\C) \cong \gr^F \C[G].
$$
The orbifold cohomology algebra $H^\hdot_{orb}([V/G],\C) =
H^\hdot_G(V,\C)^G$ is then isomorphic to $\gr^F \ZZ(G)$, where
$\ZZ(G) = \C[G]^G$ is the center of the group algebra $\C[G]$.
\end{lemma}

\proof{} Take an element $g \in G$ and consider the associated
eigenspace decomposition $V = V_1 \oplus \dots \oplus V_k$. Since $V
\cong V^*$ as a $G$-module, for every $a_i$ with $0 < a_i < r$, the
associated multiplicity $m_i$ must be the same as the multiplicity
associated to $r - a_i$. Therefore
$$
\age g = \sum_{1 \leq i \leq k}m_i\frac{a_i}{r} = \frac{1}{2}\sum_{1
\leq i \leq k,a_i \neq 0}m_i \left(\frac{a_i}{r} +
\frac{r-a_i}{r}\right) = \frac{1}{2}\codim V^g.
$$
If we have two elements $g_1,g_2 \in G$, then
$$
V^{g_1} \cap V^{g_2} \subset V^{g_1g_2}.
$$
Since for every two vector subspaces $V_1,V_2 \subset V$ we have
$\codim V_1 \cap V_2 \leq \codim V_1 + \codim V_2$, this implies
that
$$
\age g_1g_2 \leq \age g_1 + \age g_2.
$$
This proves that the filtration is indeed multiplicative. The last
claim follows.
\endproof

We note that the equality $\age g_1g_2 = \age g_1 + \age g_2$ is
achieved if two things hold: (1) the subspaces $V^{g_1}$ and
$V^{g_2}$ intersect transversaly, and (2) we have $V^{g_1g_2} =
V^{g_1} \cap V^{g_2}$ (in general, the left-hand side might be
strictly bigger). However, (1) is enough.

\begin{lemma}[{{\cite[Lemma 1.7]{GK}}}]\label{trans}
Let $V$ be a vector space equipped with a non-degenerate
skew-symmetric form $\omega$, and let $g_1$, $g_2$ be two symplectic
automorphisms of the vector space $V$. If $V^{g_1}$ and $V^{g_2}$
intersect transversally, then $V^{g_1g_2} = V^{g_1} \cap V^{g_2}$.
\end{lemma}

\proof{} We will not need really this claim, but the proof is so
simple, we could not resist reproducing it here. It is easy to see
that for any symplectic automorphism $g \in Sp(V)$, we have a
decomposition $V = V^g \oplus \Im (\id -g)$ which is orthogonal with
respect to $\omega$. Now, for any $v \in V^{g_1g_2}$ we have $g_2v =
g_1^{-1}v$, which implies
$$
(g_2-\id)v = (g_1^{-1}-\id)v.
$$
Both $g_1^{-1}$ and $g_2$ are symplectic, therefore the left-hand
side is orthogonal to $V^{g_2}$, and the right-hand side is
orthogonal to $V^{g_1^{-1}} = V^{g_1}$. Being equal, they both must
be orthogonal to $V^{g_1} + V^{g_2} = V$. We conclude that indeed
$g_2v = g_1^{-1}v = 0$.
\endproof

\section{Overview of the statements}

We will now formulate the main claims of McKay correspondence.
Firstly, there is the original cohomological McKay correspondence
conjectured in \cite{R1} and proved in \cite{B}, \cite{DL}.

\begin{theorem}[\cite{B}, \cite{DL}]\label{addi}
Let $V$ be a vector space, let $G \subset SL(V)$ be a finite
subgroup, and let $X \to V/G$ be a crepant resolution of the
quotient $V/G$. Then there exists a graded vector space isomorphism
$$
H^\hdot(X,\Q) \cong H^\hdot_{orb}([V/G],\Q).
$$
\end{theorem}

The following has also been conjectured in \cite{R1}.

\begin{corr}\label{K-th}
In the assumptions of Theorem~\ref{addi}, there exists an
isomorphism
$$
K^0_G(V) \cong K^0(X)
$$
between the $K^0$-group of $X$ and the $G$-equivariant $K^0$-group
of $V$.
\end{corr}

To derive this from Theorem~\ref{addi}, it suffices to construct a
Chern character which would identify $K_G^0(V)$ with
$H^\hdot_{orb}([V/G],\Q)$; this can be easily done by direct
inspection (it is well-known that $K^0_G(V)$ coincides with $K^0_G$
of a point, and the latter group is just the space of class
functions on the group $G$). Alternatively, one can use the very
nice general theory of Chern character and Riemann-Roch Theorem for
orbifolds developed by B. Toen in \cite{T}.

However, one could also imagine a different proof of
Corollary~\ref{K-th}, and it is for this reason that the statement was
given separately in \cite{R1}. Namely, one has the following
generalization.

\begin{conj}\label{db}
In the assumtions of Theorem~\ref{addi}, there exists an equivalence
$$
\D^b_{coh}(X) \cong D^b_G(V)
$$
between the bounded derived category of coherent sheaves on $X$ and
the bounded derived category of $G$-equivariant coherent sheaves on
$V$.
\end{conj}

So far, this remains a conjecture except in the case $\dim V=2$ (see
\cite{KV}) and the case $\dim V=3$, where a very elegant and
unexpected proof was given recently in \cite{BKR}. Unfortunately,
the methods of \cite{BKR} do not seem to apply in higher dimensions.

Finally, here is \cite[Theorem 1.2]{GK}.

\begin{theorem}[{{\cite[Theorem 1.2]{GK}}}]\label{main}
In the notation and assumptions of Theorem~\ref{addi}, assume in
addition that $V$ carries a non-degenerate skew-linear form $\omega$
preserved by $G$. Then there exists a multiplicative isomorphism
\begin{equation}\label{iso}
H^\hdot(X,\C) \cong H^\hdot_{orb}([V/G],\C).
\end{equation}
\end{theorem}

As we have noted in the Introduction, this has been proved
independently and simultaneously in \cite{V} and \cite{LS} (see also
\cite{wang}) for the case $V=\C^{2n}$, $G=S_n$ is the symmetric
group on $n$ letters, $X=\Hilb^{[n]}(\C^2)$ is the Hilbert scheme of
subschemes in $\C^2$ of dimension $0$ and length $n$. W. Wang with
co-authors were also able to treat the so-called wreath product
case: $V=\C^{2n}$, $G$ is the semidirect product of $S_n$ with the
$n$-fold product $\Gamma^n$ of a fixed finite subgroup $\Gamma
\subset SL(2,\C)$, the resolution $X$ is obtained by taking a smooth
resolution $X_0 \to \C^2/\Gamma$ and considering $X =
\Hilb^{[n]}(X_0)$.

Note that Theorem~\ref{main} in fact works for cohomology with
coefficients in $\Q$. In the symplectic case, there is an alterntive
proof of this theorem \cite{K} which also works over $\Q$. In fact,
in the Hilbert scheme case the multiplicative isomorphism
\eqref{iso} holds over $\Q$ and even over $\Z$ (possibly after one
makes certain sign changes in the definition of the orbifold
cohomology --- over $\C$, any sign changes become irrelevant). The
methods of \cite{GK}, however, very decidedly work only over $\C$.
It would be interesting to try to use the methods of \cite{K} to
check Theorem~\ref{main} over $\Q$. However, W. Wang informed me
that this might be hopeless --- his computations in the wreath
product case show that \eqref{iso} probably does not hold for
cohomology with coefficents in $\Q$. Y. Ruan \cite{Ru2} also states
his conjecture only over $\C$.

\section{Orbifold cohomology as Hochschild cohomology}

Our starting point in proving Theorem~\ref{main} is the following
observation: there is one place in mathematics where the algebra
$H^\hdot_{orb}([V/G]) \cong \gr^F \ZZ(G)$ appeared before. This is
the computation by M.S. Alvarez \cite{alv} of the Hochschild
cohomology of the so-called {\em cross-product algebra}. Let us
expain this point, starting from the definitions. We give very few
references in the next two Subsections since the subject is
well-represented in the literature; the reader is referred to, for
instanece, a comprehensive overview given in \cite{Lo}. For
everything related to qunatizations, please see \cite{Ko}.

\subsection{Recollection on Hochschild cohomology}

Let $A$ be a finitely-generated commutative algebra over $\C$, and
let $B$ be an flat associative algebra over $A$. Recall that an {\em
$B$-bimodule} $M$ is an $A$-module equipped with two commuting
$B$-module structures, one left and the other right. Equivalently,
$M$ is a left $B \otimes_A B^{op}$-module, where $B^{op}$ denotes
the opposite algebra. All $B$-bimodules form an abelian
category. The algebra $B$ itself is tautologically a
$B$-bimodule. For any $B$-bimodule $M$, one defines the {\em
Hochschild cohomology groups} $HH^\hdot_A(B,M)$ with coefficients in
$M$ by
$$
HH^\hdot_A(B,M) = \Ext^\hdot_{B \otimes_A B^{op}}(B,M).
$$
The Hochschild cohomology groups $HH^\hdot_A(B,B)$ with coefficients
in the tautological module $B$ are denoted simply $HH^\hdot_A(B)$
and called the {\em Hochschild cohomology groups of the algebra
$B$}.

This notion immediately generalizes to the case of sheaves of
algebras. Namely, if $X$ is a topological space equipped with a
sheaf $\calo$ of flat associative $A$-algebras, then by definition
$$
HH^\hdot_A(X) = \Ext^\hdot_{X \times
X^{op}}(\calo_\Delta,\calo_\Delta),
$$
where $\Ext$-groups are taken in the category of sheaves of $\calo
\boxtimes \calo^{op}$-modules on $X \times X$, and $\calo_\Delta$ is
the tautological sheaf supported on the diagonal $\Delta \subset X
\times X$. In particular, in this way one defines Hochschild
cohomology of flat scheme over $\Spec A$. Even further, one can
define Hochschild cohomology of an arbitraty abelian $A$-linear
category, and the resulting groups are {\em dervied
Morita-invariant} -- this means that if two $A$-linear abelian
categories have equivalent derived categories, then their Hochschild
cohomology groups are naturally isomorphic.

Hochschild cohomology $HH^\hdot_A(B)$ carries two natural algebraic
structures. The first one is the multiplication given by the Yoneda
multiplication in the $\Ext$-groups. The second one is a Lie bracket
$\{-,-\}$ of degree $-1$ called {\em the Gertsenhaber bracket}. It
is not at all obvious from the definition, but it exists
nevertheless; moreover, the bracket and the multiplication together
form a so-called {\em Gerstenhaber algebra}. We do not list the
axioms here because we will not need them. We only remark that one
of them says that the Yoneda multiplication is commutative.

In the case when the algebra $B$ is commutative, finitely generated
and smooth over $A$, the Hochschild cohomology $HH^\hdot_A(B)$ has
been computed in the classic paper \cite{HKR}. The answer is
$$
HH^\hdot_A(B) = \Lambda^\hdot_A\T(B/A),
$$
the exterior algebra generated by the module $\T(B/A)$ of $A$-linear
derivations $\delta:B \to B$. In other words, Hochshchild cohomology
classes are the same as polyvector fields on $B/A$. The
multiplication is just the exterior algebra multiplication, and the
Gerstenhaber bracket is the so-called {\em Schouten bracket} of
polyvector fields. This result gives rise to a useful modification
of the Hochschild cohomology known as {\em Poisson
cohomology}. Recall that a {\em Poisson structre} on $B$ over $A$ is
given by an $A$-linear Lie bracket $\{-,-\}:B \otimes_A B \to B$
which is a derivation with respect to the multiplication in
$B$. Equivalently, a Poisson structure is given by bivector field
$\Theta \in \Lambda^2_A\T(B/A)$ such that $\{\Theta,\Theta\}=0$ with
respect to the Schouten bracket (the Poisson bracket is then
$\{f,g\} = \langle \Theta, df \wedge dg \rangle$). In the Hochschild
cohomology language, $\Theta$ becomes a class in $HH^2_A(B)$ such
that $\{\Theta,\Theta\}=0$ with respect to the Gerstenhaber
bracket. Given such a class, one can introduce a map
$\delta_\Theta:HH^\hdot_A(B) \to HH^{\hdot+1}_A(B)$ by setting
$$
\delta_\Theta(a) = \{\Theta,a\}.
$$ 
This has been done by J.-L. Brylinski \cite{Br}. It follows from the
axioms of the Gerstenhaber algebra, -- or equivalently, from the
properties of the Schouten bracket, -- that $\{\Theta,\Theta\}=0$
implies $\delta_\Theta^2=0$, and that $\delta_\Theta$ is a
derivation with respect to the multiplication in
$HH^\hdot_A(B)$. Therefore we obtain a differential-graded algebra
$\langle HH^\hdot_A(B), \delta_\Theta \rangle$. Its cohomology
groups are denoted by $HP^\hdot_A(B)$ and called {\em the Poisson
cohomology groups} of the Poisson algebra $B/A$.

We note that this definition makes sense for any $B$, not only for a
smooth and commutative one. If $B$ is not commutative, it no longer
makes sense to speak of a Poisson structure on $B$. Nevertheless,
for any $\Theta \in HH^2_A(B)$ with $\{\Theta,\Theta\} = 0$ we do
have the differential-graded algebra $\langle HH^\hdot_A(B),
\delta_\Theta \rangle$. We will denote its cohomology groups by
$HH^\hdot_\Theta(B)$ and call them the {\em twisted Hochschild
cohomology groups}. Note that in the case when $B$ is commutative
but not smooth, the notion of a Poisson structure on $B$ is still
well-defined. There also exists a natural notion of Poisson
cohomology groups $HP^\hdot_A(B)$; however, this notion is more
complicated, and the groups themselves in general do {\em not}
coincide with the twisted Hochschild cohomology groups
$HH^\hdot_\Theta(B)$ (a large part of \cite[Appendix]{GK} is devoted
to a detailed investigation of these phenomena).

\subsection{Recollection on quantizations.}

One situation where twisted Hoch\-schild cohomology occurs naturally
is the following one. Let $B$ be an associative algebra over
$\C$. Take $A = \C[h]$, the algebra of polynomials in one variable
which we denote by $h$, and let $B_h$ be a flat associative
$A$-algebra such that $B_h/h \cong B$ -- in other words, $B_h$ is a
one-parameter deformation of the algebra $B$. Then the deformation
theory associates a Hochschild cohomology class $\Theta \in
HH^2_\C(B)$ to the family $B_h$. One can compute the Hochschild
cohomology $HH^\hdot_A(B_h)$ by using the spectral sequence
associated to the $h$-adic filtration. The term $E^1$ of this
spectral sequence coincides with $HH^\hdot_\C(B) \otimes_\C A$, and
the differential is given by $\delta_\Theta$ (in particular, we have
$\{\Theta,\Theta\} = 0$). Thus the term $E^2$ of the spectral
sequence coincides (for dimensional reasons, only modulo
$h$-torsion) with the twisted Hochschild cohomology groups
$HH^\hdot_\Theta(B) \otimes_\C A$.

If the algebra $B$ is commutative, then the class $\Theta$ induces
a Poisson structure on $B$. The Poisson bracket on $B$ is given
by
$$
\{f,g\} = \frac{1}{h}\left(\wt{f}\wt{g} - \wt{g}\wt{f}\right) \mod
h^2,
$$
where $\wt{f}$, $\wt{g}$ are arbitrary elements in $B_h$ with
$\wt{f}=f \mod h$, $\wt{g}=g \mod h$. The algebra $B_h$ is then said
to be a {\em quantization} of the Poisson algebra $B$.

We will need a particular case of these general constructions,
namely, the case when the commutative algebra $B$ is smooth and {\em
symplectic} -- in other words, the skew-symmetric $A$-bilinear
pairing $\Omega^1(B) \otimes_B \Omega^1(B) \to B$ on the cotangent
module $\Omega^1(B)$ given by Poisson bivector $\Theta \in
\Lambda^2\T(B)$ is non-degenerate. In this case, $\Theta$ induces an
isomorphism $\Omega^1(B) \to \T(B)$, and one can identify the
modules $\Lambda^\hdot\T(B)$ of polyvector fields on $B$ with the
modules $\Omega^\hdot(B)$ of differential forms on $\Spec B$. Under
this identification, -- and this one of the important results of
\cite{Br}, -- the Poisson differential
$\delta_\Theta:\Lambda^\hdot\T(B) \to \Lambda^{\hdot+1}\T(B)$
becomes the de Rham differential $d:\Omega^\hdot(B) \to
\Omega^{\hdot+1}(B)$. Therefore the Poisson cohomology $HP^\hdot(B)$
concides with the de Rham cohomology $H^\hdot_{DR}(\Spec B)$.

The simplest example of a smooth symplectic algebra is the algebra
$B = S^\hdot(V)$ of polynomials on a symplectic vector space $V$
with symplectic form $\omega$. Then $HP^k(B) = H^k_{DR}(V)$ is $\C$
for $k=0$ and $0$ otherwise. The algebra $S^\hdot(V)$ has a standard
qunatization $W_h$ known as the {\em Weyl algebra}; it is the
associative algebra over $A = \C[h]$ generated by $V$ modulo the
relations
$$
vw-wv = h\omega(v,w)
$$
for all $v,w \in V$. Specializing to $h=1$, we get the Weyl algebra
$W = W_h/(h-1)$. If $\dim V =2n$, then $W$ is non-canonically
isomorphic to the algebra of differentail operators on the affine
space $\A^n$. Since $HP^k(B) = 0$ for $k \geq 1$, the $h$-adic
spectral sequence for $HH^\hdot_A(W_h)$ degenerates. We deduce that
$HH^0(W) = \C$ and $HH^k(W) = 0$ for $k \geq 1$.

\subsection{Computation for the smash-product algebra.}

For any associative algebra $A$ equipped with an action of a finite
group $G$, introduce the {\em smash-product algebra} $A \# G$ in the
following way: $A \# G = A \otimes_\C \C[G]$ as a vector space, and
we have
$$
(a_1 \cdot g_1)(a_2 \cdot g_2) = a_1(a_2^g) \cdot g_1g_2
$$
for all $a_1,a_2 \in A$ and $g_1,g_2 \in G$. In particular, given a
symplectic vector space $V$ and a finite subgroup $G \subset Sp(V)$,
we can form the smash-product $W \# G$, where $W$ is the Weyl
algebra ssociated to $V$. In \cite{alv}, M.S. Alvarez has proved the
following.

\begin{theorem}[\cite{alv}]\label{alva}
The Hochschild cohomology algebra $HH^\hdot(W \# G)$ is isomorphic
$$
HH^\hdot(W \# G) \cong \gr^F \ZZ(G)
$$
to the associated graded $\gr^F \ZZ(G)$ of the center $\ZZ(G)$ of
the group algebra $\C[G]$ with respect to the filtration introduces
in Lemma~\ref{filtr}.
\end{theorem}

This Theorem is our main motivation for Theorem~\ref{main}. Its
proof, of which we will only present a sketch, is based on an
earlier result of \cite{alev}, where $HH^\hdot(W \# G)$ has been
computed as a graded vector space. For simplicity, we will first do
this compuation with $W$ replaced by the commutative polynomial
algebra $B=S^\hdot(V)$.

\begin{lemma}
Consider the smash-product algebra $B \# G$. Then there exist an
isomorphism
\begin{equation}\label{alv.eq}
HH^\hdot(B \# G) = \left( \bigoplus_{g \in
G}\Omega^\hdot(V^g)[\codim V^g]\right)^G,
\end{equation}
where for any $g \in G$ we denote by $V^g \subset V$ the subspace of
$g$-invariant vectors in $V$, and $\Omega^k(V^g)$ is the space of
differential forms of degree $k$ on the vector space $V^g$.
\end{lemma}

\proof[Sketch of the proof.] Since the group $G$ is finite, and our
base field $\C$ is of characteristic $0$, the group $G$ has no
cohomology. Therefore we have
\begin{align*}
HH^\hdot(B \# G) &= \Ext^\hdot_{B\#G-\bimod}(B\#G,B\#G)\\ &=
\left(\Ext^\hdot_{B-\bimod}(B,B\#G)\right)^G\\ &=
\left(HH^\hdot(B,B\#G)\right)^G
\end{align*}
The $B$-bimodule $B\#G$ can be decomposed
$$
B\#G = \bigoplus_{g \in G}B^g,
$$
where $B^g$ is $B$ with the standard right $B$-module structure, and
with the left $B$-module structure twisted by the action of $g \in
G$. 

Geometrically, finitely-generated $B$-bimodules are coherent sheaves
on $V \times V = \Spec(B \otimes_\C B^{op})$; then $B^g$ is the just
the structure sheaf of the graph $\graph g \subset V \times V$ of
the map $g:V \to V$. In particular, $B = B^\id$ itself is the
structure sheaf of the diagonal $\Delta \subset V \times V$. To
prove the Lemma, it remains to show that
$$
HH^\hdot(B,B^g) = \Omega^\hdot[V^g][\codim V^g].
$$
This is immediate. Indeed, $V^g = \Delta \cap \graph g$, and $V^g$
is a symplectic vector space. One computes $HH^\hdot(B,B^g)$ by the
local-to-global spectral sequence on $V \times V$, and then
identifies polyvector fields on $V^g$ with the differential forms by
means of the symplectic structure.
\endproof

To deduce Theorem~\ref{alva}, one first uses the $h$-adic spectral
sequence to compute $HH^\hdot(W\#G)$ as a vector space. At the term
$E^2$, we have the twisted Hochschild cohomology groups
$HH^\hdot_\Theta(B^g)$. One observes immediately that the
differential $\delta_\Theta$ coincides with the de Rham differential
on $V^g$, so that $HH^\hdot_\Theta(B,B^g) = H^\hdot_{DR}(V^g)[\codim
V^g]$. The spectral sequence degenerates at $E^2$ (for instance,
because all the non-trivial classes in $E^2$ have even degrees). We deduce
$$
HH^\hdot(W\#G) = \left(\bigoplus_{g \in G}\C[\codim V^g]\right)^G,
$$
which is exactly \eqref{alv.eq}. The multiplicative structure is not
easy to see from the spectral seqeunce. To see it, it is simpler to
represent the classes generating $HH^\hdot(W\#G)$ by explicit
differential forms $\omega_g$ (this can be made even more precise
and explicit if one uses the Koszul complex for the vector space $V$
to resolve $B$ as a $B$-module and compute $HH^\hdot(B,B^g)$). Then
one immediately checks that for every $g.h \in G$, we have
$\omega_g\omega_h=\omega_{gh}$ if $V^g, V^h \subset V$ intersect
transversally, and $\omega_g\omega_h = 0$ otherwise. By
Lemma~\ref{trans}, this is exactly the multiplicative structure in
$\gr^F\ZZ(G)$. Note that in fact the same algebra appears as the
twisted Hochschild cohomology algebra
$HH^\hdot_\Theta(S^\hdot(V)\#G)$.

\section{The proofs.}

\subsection{A proof.}

We will now present what ought to be the proof of
Theorem~\ref{main}. It will only take half a page.

Let $V$ be a symplectic vector space, let $G \subset Sp(V)$ be a
finite subgroup, and let $X \to V/G$ be a crepant resolution of
singularities of the quotient $V/G$. Since the symplectic form on
$V$ is $G$-invariant, it descends to a symplectic form on the smooth
part of the quotient $V/G$. This form then extends to a closed
non-degenerate $2$-form on the resolution $X$ (see
e.g. \cite{K1}). Therefore the variety $X$ is equipped with a
natural Poisson structure $\Theta$, and we can compute its de Rham
cohomology by comparison with the Poisson cohomology. We then have
the isomorphisms
\begin{equation}\label{first}
H^\hdot_{DR}(X) \cong HP^\hdot(X) \cong HH^\hdot_\Theta(X).
\end{equation}
Now, {\em assume that Conjecture~\ref{db} is true}. Note that the
category of $G$-equivariant coherent sheaves on $V$ is equivalent to
the category of finitely-generated modules over the smash-product
algebra $B\#G=S^\hdot(V)\#G$. Since the Hochschild cohomology is
dervied Morita-invariant, we have
\begin{equation}\label{second}
H^\hdot_\Theta(X) \cong HH^\hdot_\Theta(B\#G).
\end{equation}
This finishes the proof: indeed, by Theorem~\ref{alva} (or rather,
by its proof) we have $HH^\hdot_\theta(B\#G) \cong
H^\hdot_{orb}([V/G],\C)$.

\subsection{A debunking.} Having given a proof of
Theorem~\ref{main}, we will now demolish it.

The most obvious problem with the proof is its reliance on
Conjecture~\ref{db}. However, this is not so serious. By a stroke of
luck, there exists an approach to this conjecture, which is now
under investigation (and the provisonary reference for this is
\cite{BK}). The author is reasonably certain that in the nearest
future Conjecture~\ref{db} will be proved; hopefully this will
happen by the time the present volume is out of print.

Unfortunately, our fake proof of Theorem~\ref{main} also has many
other flaws, and we will now enumerate those.
\begin{enumerate}
\item We have used isomorphisms \eqref{first} in spite of the fact
that the manifold $X$ is not affine.

This causes two problems. Firstly, the terms in \eqref{first} might
not be defined for non-affine $X$. As we have explained, the
definition of the {\em algebra} $HH^\hdot(X)$ does not require $X$
to be affine. To define the twisted version, however, we also need
the Gerstenhaber bracket. This is possible to define, too:
essentially, the Hochschild cohomology algebra $HH^\hdot(X)$ can be
obtained as the hypercohomology of a certain complex of sheaves on
$X$, this complex carries the structure of a Gerstenhaber algebra,
the standard homotopy techniques as in e.g. \cite{HS} give a
homotopy Gerstenhaber algebra on $HH^\hdot(X)$, and this is
equivalent to a usual Gerstenhaber algebra structure by Kontsevich
formality (for this point, we would recommend to consult \cite{HT}
and references therein). This is really quite roundabout, one
certainly would prefer a more direct approach. This is currently the
topic of active research. At least one approach has been suggested
recently in \cite{R}, \cite{ke2}.

Secondly, even assuming that all the algebras in \eqref{first} are
well-defined, they might not be isomorphic when $X$ is not
affine. And this is exactly what happens, unfortunately: the
Hochschild-Kostant-Rosenberg isomorphism $HH^\hdot(X) \cong
H^\hdot(X,\Lambda^\hdot\T(X))$ holds for non-affine varieties, but
{\em it is not compatible with the Gerstenhaber algebra structure}.

It is probably true that this isomorphism can be {\em corrected} in
a certain precise way so that it becomes multiplicative -- this has
been claimed without proof in the last part of \cite{Ko}, and has
been the subject of much research (I do not feel competent enough to
provide an exact reference, except for saying that for {\em compact}
$X$ the problem has been definitely completely solved by
A. Caldararu \cite{C2}). However, whether this isomorphism is also
compatible with the Gerstenhaber bracket is anyone's guess.

\item We have used the derived Morita-invariance of Hochschild
homology in \eqref{second} without any justification.

The derived Morita-invariance property for Hochschild {\em homology}
has been proved several years ago by B. Keller, see \cite{ke1}. The
invariance of cohomology is much simpler. It is so simple that
again, I have not been able to track a reference. It might be that
no one was diligent enough to write the proof down. For compact $X$
(and derived equivalences of Fourier-Mukai type) everything has been
recently written down very carefully by A. Caldararu \cite{C1} (as
the foundation for his main results). But again, the question of the
Gerstenhaber bracket seems to be open. Quite recently (two weeks
ago, at the time of the writing) B. Keller has published a preprint
\cite{ke2} where this question might be solved; unfortunately, my
lack of real expertise in the subject forces me to reserve judgement
at this time.

\item We have tacitly assummed that the twisting cocycle $\Theta$ on
$B\#G$ obtained from $X$ is the same as the standard cocycle
obtained by descent from the symplectic form on $V$.

This should be fairly easy to check, {\em once the Morita-invariance
has been solidly established}. Unfortunately, we cannot even start
this checking before we know exactly how the isomorphism
\eqref{second} works.
\end{enumerate}

\subsection{How the real proof works.} We will now describe, very
briefly and sketchily, the main stages of the real proof of
Theorem~\ref{main} given in \cite{GK} (a somewhat more detailed
description can be found in \cite[Introduction]{GK}). The main idea
is the following: all the difficulties appear because we need to
work with a non-affine variety $X$. Things would be much simpler if
one could {\em deform} everything so that the manifold in question
becomes affine.

\medskip

On the resolution side of the picture, a reasonably good deformation
theory for smooth non-compact symplectic varieties has been
developed in \cite{kave}. In parituclar, it can be applied to a
crepant resolution $X \to V/G$. The result is a smooth family
$\wt{X}/B$ parametrized by the affine space $B=H^2_{DR}(X)$. The
approach in \cite{kave} only gives a formal deformation; however,
there is a natural $\C^*$-action on $V/G$ and on $X$ which allows
one to spread out the deformation to the whole $B$. The de Rham
cohomology $H^\hdot_{DR}(\wt{X}_b)$ of the fiber $\wt{X}_b$ is the
same for every point $b \in B$. The deformation $\wt{X}$ also
induces a deformation $\wt{Y}/B$ of the quotient $Y = V/G$, which is
affine. The map $X \to V/G$ extends to a map $\wt{X} \to \wt{Y}$,
and the extended map is {\em one-to-one over a generic point $b \in
B$}. Consequently, for generic $b \in B$ the variety
$\wt{X}_b=\wt{Y}_b$ is simultaneously smooth (being part of
$\wt{X}/B$) and affine (being part of $\wt{Y}/B$).

\medskip

On the smash-product side, a very nice deformation $H_{t,c}$ of the
smash-product algebra $S^\hdot(V)\#G$ has been constructed and
studied in some detail in \cite{EG}. The parameter space for this
deformation is the product of an affine line with coordinate $t$ and
an affine space $C$ of $G$-invariant functions on the set $S \subset
G$ of elements $g \in G$ of age $1$ (these elements are also known
as {\em symplectic reflections}). The deformation in the
$t$-direction is entirely non-commutative; in particular, it
incorporates the Weyl algebra deformation $W_h\#G$. The deformation
$H_{0,c}$, while still non-commutative, is closer to the commutative
world -- namely, the center $Z_c$ of the algebra $H_{0,c}$ gives a
flat deformation of the center $Z_0 \subset S^\hdot(V)\#G$. The
latter is obviously just the subalgebra of $G$-invariant polynomials
in $S^\hdot(V)$ -- in other words, the algebra of algebraic
functions on the quotient $V/G$.

Moreover, the algebras $Z_c$ carry a natural Poisson structure
(which is essentially obtained from the additional deformation in
the $t$-direction), and all of the steps of our fake proof have been
solidly established in \cite{EG} {\em for the affine variety $\Spec
Z_c$}. In particular, Conjecture~\ref{db} becomes a proposition
which claims that the algebra $H_{0,c}$ is Morita-equivalent to its
center $Z_c$. Some of the claims -- including this
Morita-equivalence -- require one to know {\em a priori} that $\Spec
Z_c$ is smooth. If this is known, then the de Rham cohomology
$H^\hdot_{DR}(\Spec Z_c)$ has been computed in \cite{EG}, and the
answer is $H^\hdot_{orb}([V/G],\C)$.

\medskip

These two deformation were the starting point for \cite{GK}. We note
that by the additive McKay correspondence, the base spaces
$B=H^2_{DR}(X)$ and $C$ are canonically identified. One might hope
to identify the deformations $\Spec Z_c$ and $\wt{Y}/B$. Then one
deduces that $\Spec Z_c$ is smooth for generic $c \in C$, and
applies \cite{EG} to compute its cohomology.

Unfortunately, we were not able to prove that such an identification
exists, and we had to settle for less. Namely, we prove that
$\wt{Y}/B$ gives a {\em versal} deformation of the quotient $Y=V/G$
in the class of affine Poisson schemes. In other words, a fiber of
any Poisson deformation also occurs as a fiber of the deformation
$Z_c$, but possibly in a non-unique way.  In fact, there exists a
universal deformation, too, and its base $S$ is also a smooth affine
space of dimension $\dim S = \dim B = \dim C$. However, the
classifying maps $B \to S$, $C \to S$ of the deformations $\wt{Y}/B$
and $Z_c/C$ are {\em ramified covers}.

The model for this situation is the case when $\dim V = 2$; then the
subgroups $G \subset Sp(V) = SL(2,\C)$ are classified by simple Lie
algebras $\g$ of types $A$, $D$, $E$, both the space $B$ and the
space $C$ are naturally identified with a Cartan subalgebra $\h
\subset \g$, and the space $S$ is the quotient $\h/W$ of the Cartan
algebra $\h$ by the Weyl group $W$. We expect that the picture in
higher dimension is similar, -- in particular, there should exist a
natural identification $B = C$. But we are not able to make a
precise conjecture at this time.

Be that as it may, what we can prove, eventually, is that every
fiber $\Spec Z_c$ of the deformation $Z_c$ occurs as a fiber
$\wt{Y}_b$ of the deformation $\wt{Y}$, and $b$ is generic when $c$
is generic. This shows that $\Spec Z_c$ is smooth for generic $c$
and that 
$$
H^\hdot_{DR}(\Spec Z_c) \cong H^\hdot_{DR}(\wt{Y}_b) \cong
H^\hdot_{DR}(\wt{X}_b) \cong H^\hdot_{DR}(X).
$$
We then invoke \cite{EG} to compute the left-hand side and prove
Theorem~\ref{main}.

The only new ingredient in \cite{GK} as compared to \cite{kave} and
\cite{EG} is a deformation theory for Poisson algebras which is
developed far enough to prove the versality of the deformation
$\wt{Y}/B$.

\begin{remark} To conclude the paper, we would like to note that the idea to
relate orbifold cohomology to something in the Hochschild cohomology
world has been introduced by V. Baranovsky \cite{ba}. He works in
full generality -- an arbitrary Calabi-Yau quotient $V/G$ -- but he
used cyclic (or Hochschild) homology, not Hochschild cohomology. It
is certainly a more natural approach. Unfortunately, the homology
groups do not carry a multiplicative structure. To try to rectify
this, one may attempt to identify Hochschild homology and Hochschild
cohomology by cupping a cohomology class with the volume
form. However, {\em this is not what we do}. Our proof also contain
a hidden identification $HH_\idot \cong HH^\hdot$, but the
identification uses the symplectic form. On the level of polyvector
fields, we use the isomorphism $\T(X) \cong \Omega^1(X)$ -- not the
isomorphism $\T(X) \cong \Omega^{2\dim X -1}(X)$ obtained by cupping
with the volume form.
\end{remark}

\bigskip

\noindent
{\sc Steklov Math Institute\\
Moscow, USSR}

\bigskip

\noindent
{\em E-mail address\/}: {\tt kaledin@mccme.ru}


\begin{thebibliography}{DWL}

\bibitem[AFLS]{alev} J. Alev, M.A. Farinati, T. Lambre, and
A.L. Solotar, {\it Homologie des invariants d'une alg\`ebre de Weyl
sous l'action d'un groupe fini}, J. of Algebra {\bf 232} (2000),
564--577.

\bibitem[Al]{alv} M.S. Alvarez, {\em Algebra structure on the
 Hochschild cohomology of the ring of invariants of a Weyl algebra
 under a finite group}, J. Algebra {\bf 248} (2002), 291--306.

\bibitem[Bar]{ba} V. Baranovsky, {\em Orbifold Cohomology as
Periodic Cyclic Homology}, math.AG/0206256.

\bibitem[Bat]{B} V. Batyrev, {\em Non-Archimedian integrals and
stringy Euler numbers of log-terminal pairs}, Jour. of European
Math. Soc. {\bf 1} (1999), 5--33.

\bibitem[BK]{BK} R. Bezrukavnikov and D. Kaledin, {\em McKay
equivalence for symplectic resolutions of singularities}, in
preparation.

\bibitem[BKR]{BKR} T. Bridgeland, A. King, and M. Reid, {\em The
McKay correspondence as an equivalence of derived categories,}
Jour. AMS {\bf 14} (2001), 535--554.

\bibitem[Br]{Br} J.-L. Brylinski, {\it A differential complex for
Poisson manifolds}, J. Diff. Geom. {\bf 28} (1988), 93--114.

\bibitem[C1]{C1} A. Caldararu, {\em The Mukai pairing, I: the
Hochschild structure}, math.AG/0308079.

\bibitem[C2]{C2} A. Caldararu, {\em The Mukai pairing, II: the
Hochschild-Kostant-Rosenberg isomorphism}, math.AG/0308080.

\bibitem[CR]{CR} W. Chen and Y. Ruan, {\em A new cohomology theory
for orbifold}, math.AG/0004129.

\bibitem[DL]{DL} J. Denef and F. Loeser, {\em Motivic integration,
quotient singularities and the McKay correspondence}, Comp.
Math. {\bf 131} (2002), 267--290.

\bibitem[EG]{EG} P. Etingof and V. Ginzburg, {\em Symplectic
reflection algebras, Calogero-Moser space, and deformed
Harish-Chandra homomorphism}, Invent. Math. {\bf 147} (2002),
243--348.

\bibitem[FG]{FG} B. Fantechi and L. G\"ottsche, {\em Orbifold
cohomology for global quotients}, Duke Math. J. {\bf 117} (2003),
197--227.

\bibitem[GK]{GK} V. Ginzburg and D. Kaledin, {\em Poisson
deformations of symplectic quotient singularities}, math.AG/0212279.

\bibitem[H]{HT} V. Hinich, {\em Tamarkin's proof of Kontsevich
formality theorem},\\ math.QA/0003052.

\bibitem[HS]{HS} V. Hinich and V. Schechtman, {em Deformation theory
and Lie algebra homology}, alg-geom/9405013.

\bibitem[HKR]{HKR} G. Hochschild, B. Kostant, and A. Rosenberg, {\em
   Differential forms on regular affine algebras}, Trans. AMS {\bf
   102} (1962), 383--408.

\bibitem[Ka1]{K1} D. Kaledin, {\em Dynkin diagrams and crepant
  resolutions of quotient singularities}, math.AG/9903157.

\bibitem[Ka2]{K} D. Kaledin, {\em McKay correspondence for
symplectic quotient singularities}, Invent. Math. {\bf 148} (2002),
151--175.

\bibitem[KaVe]{kave} D. Kaledin and M. Verbitsky, {\em Period map
   for non-compact holomorphically symplectic manifolds}, GAFA {\bf
   12} (2002), 1265--1295.

\bibitem[KaVa]{KV} M. Kapranov and E. Vasserot, {\em Kleinian
    singularities, derived categories and Hall algebras}, Math. Ann.
    {\bf 316} (2000), 565--576.

\bibitem[Ke1]{ke1} B. Keller, {\em On the cyclic homology of exact
categories}, J. Pure Appl. Algebra {\bf 136} (1999), 1--56.

\bibitem[Ke2]{ke2} B. Keller, {\em Hochschild cohomology and derived
  Picard groups},\\ math.KT/0310221.

\bibitem[Ko]{Ko} M. Kontsevich, {\em Deformation quantization of
  Poisson manifolds, I}, q-alg/9709040.

\bibitem[LS]{LS} M. Lehn and Ch. Sorger, {\em Symmetric groups and
the cup product on the co\-ho\-mo\-logy of Hilbert schemes}, Duke
Math. J. {\bf 110} (2001), 345--357.

\bibitem[LQW]{wang} W. Li, Z. Qin, and W. Wang, {\em Vertex algebras
and the cohomology ring structure of Hilbert schemes of points on
surfaces}, Math. Ann. {\bf 324} (2002), 105--133.

\bibitem[L]{Lo} J.-L. Loday, {\em Cyclic homology}, Gr\"undlehren
der Mathematischen Wissenschaften {\bf 301}, Springer-Verlag,
Berlin, 1992.

\bibitem[McK]{McK} J. McKay, {\em Graphs, singularities and finite
groups}, in {\em The Santa Cruz Conference on Finite Groups},
Proc. Symp. Pure Math. {\bf 37} (1980), 183--186.

\bibitem[R1]{R1} M. Reid, {\em McKay correspondence},
alg-geom/9702016 v3, 1997.

\bibitem[R2]{R2} M. Reid, {\em La correspondance de McKay},
S\'eminaire Bourbaki, Vol. 1999/2000. Ast\'erisque No.  {\bf 276}
(2002), 53--72.

\bibitem[RZ]{R} R. Rouquier and A. Zimmerman, {\em Picard groups for
derived module categories}, Proc. Lond.Math.Soc. {\bf 87} (2003),
197--225.

\bibitem[Ru1]{Ru1} Y. Ruan, {\em Stringy Geometry and Topology of
Orbifolds}, Symposium in Honor of C. H.  Clemens, 187--233,
Contemp. Math., {\bf 312}, Amer. Math. Soc., Providence, RI, 2002.

\bibitem[Ru2]{Ru2} Y. Ruan, {\em Stringy Orbifolds}, Orbifolds in
mathematics and physics (Madison, WI, 2001), 259--299,
Contemp. Math., {\bf 310}, Amer. Math. Soc., Providence, RI, 2002.

\bibitem[T]{T} B. Toen, {\em Th\'eor\`emes de Riemann-Roch pour les
champs de Deligne-Mumford}, $K$-Theory {\bf 18} (1999), 33--76.

\bibitem[V]{V} E. Vasserot, {\em Sur l'anneau de cohomologie du
sch\`ema de Hilbert de $\mathbf C\sp 2$}, C.R. Acad. Sci. Paris,
{\bf 332} (2001), 7--12.

\end{thebibliography}
\end{document}